\documentclass[11pt]{article}
\usepackage[T1]{fontenc}
\usepackage{lmodern}
\usepackage[margin=0.88in]{geometry}
\usepackage{microtype}
\usepackage{amsmath,amssymb,amsthm}
\usepackage{booktabs,array}
\usepackage{xcolor}
\usepackage{float}
\usepackage{placeins}
\floatstyle{ruled}
\newfloat{algorithm}{htbp}{loa}
\floatname{algorithm}{Algorithm}
\newcounter{algline}
\newcount\algdepth
\newenvironment{algorithmic}[1][1]{%
  \small\setcounter{algline}{0}\algdepth=0
  \setlength{\parskip}{1.4pt}\setlength{\parindent}{0pt}
}{\par}
\newcommand{\State}{\par\stepcounter{algline}\noindent
  \makebox[1.6em][r]{\scriptsize\thealgline:}\hspace{.5em}%
  \hspace*{\dimexpr\algdepth em\relax}\ignorespaces}
\newcommand{\Require}{\par\noindent\textbf{Input:}\ }
\newcommand{\If}[1]{\State\textbf{if} #1 \textbf{then}\advance\algdepth by1\relax}
\newcommand{\EndIf}{\advance\algdepth by-1\relax\State\textbf{end if}}
\newcommand{\For}[1]{\State\textbf{for} #1 \textbf{do}\advance\algdepth by1\relax}
\newcommand{\While}[1]{\State\textbf{while} #1 \textbf{do}\advance\algdepth by1\relax}
\newcommand{\EndWhile}{\advance\algdepth by-1\relax\State\textbf{end while}}
\newcommand{\EndFor}{\advance\algdepth by-1\relax\State\textbf{end for}}
\newcommand{\Function}[2]{\State\textbf{function} \textsc{#1}(#2)\advance\algdepth by1\relax}
\newcommand{\EndFunction}{\advance\algdepth by-1\relax\State\textbf{end function}}
\newcommand{\Return}{\textbf{return}\ }
\newcommand{\Call}[2]{\textsc{#1}(#2)}

\usepackage{tikz}
\usepackage{needspace}
\usepackage[hidelinks]{hyperref}
\usepackage{authblk}
\hypersetup{pdftitle={Exact Values and Bounds for Four Color Quadratic Polynomial van der Waerden Numbers},pdfauthor={Ishan Kumthekar, William Gasarch, William Regli, Ian Whitehouse}}

\AtBeginDocument{%
\setlength{\abovedisplayskip}{6pt plus 2pt minus 2pt}
\setlength{\belowdisplayskip}{6pt plus 2pt minus 2pt}
\setlength{\abovedisplayshortskip}{4pt plus 2pt minus 2pt}
\setlength{\belowdisplayshortskip}{4pt plus 2pt minus 2pt}}
\newtheorem{theorem}{Theorem}[section]
\newtheorem{lemma}[theorem]{Lemma}
\newtheorem{proposition}[theorem]{Proposition}
\theoremstyle{definition}
\newtheorem{definition}[theorem]{Definition}

\theoremstyle{remark}

\newcommand{\Z}{\mathbb Z}
\newcommand{\abs}[1]{\left|#1\right|}
\newcommand{\Dist}[1]{\mathcal D\!\left(#1\right)}
\newcommand{\PVDW}[2]{W\!\left(#1;#2\right)}
\theoremstyle{plain}
\newtheorem{corollary}[theorem]{Corollary}
\colorlet{witnessred}{red!82!black}
\colorlet{witnessgreen}{green!45!black}
\colorlet{witnessblue}{blue!82!black}
\colorlet{witnessyellow}{yellow!58!black}
\tikzset{witness point/.style={circle,draw=white,line width=.65pt,minimum size=5.4mm,inner sep=0pt},witness edge/.style={draw=black!38,line width=.45pt}}
\newcommand{\Rcol}{\textcolor{red!82!black}{\textbf{R}}}
\newcommand{\Gcol}{\textcolor{green!45!black}{\textbf{G}}}
\newcommand{\Bcol}{\textcolor{blue!82!black}{\textbf{B}}}
\newcommand{\Ycol}{\textcolor{yellow!58!black}{\textbf{Y}}}
\newcommand{\witnessbox}[1]{%
  \begin{center}%
  \setlength{\fboxsep}{3pt}%
  \fbox{\begin{tabular}{@{}c@{}}#1\end{tabular}}%
  \end{center}%
}

\title{Exact Values and Bounds for Four Color Quadratic Polynomial van der Waerden Numbers}
\author[1]{Ishan Kumthekar}
\author[2]{William Gasarch}
\author[2]{William Regli}
\author[2]{Ian Whitehouse}

\affil[1]{Department of Mathematics, University of Florida,
Gainesville, FL 32611, USA\\
\texttt{ishankumthekar@ufl.edu}}

\affil[2]{Department of Computer Science, University of Maryland,
College Park, MD 20742, USA\\
\texttt{gasarch@umd.edu, regli@umd.edu, ianjw@umd.edu}}

\date{September 5, 2026}
\begin{document}
\maketitle
\vspace{-1.2em}
\begin{abstract}
For a nonzero polynomial $P(x)\in\mathbb{Z}[x]$ with $P(0)=0$, let $W(P(x);4)$ denote the least integer $N$ such that every four coloring of $[N]$ contains distinct integers $u,v$ of the same color with $v-u=P(d)$ for some integer $d$. The best known general upper bounds for polynomial van der Waerden numbers are primitive recursive but remain enormous compared with the few known exact values.

In this paper we determine the exact four color polynomial van der Waerden numbers for
\[
\begin{gathered}
x^2+x,\quad 2x^2+x,\quad 3x^2+x,\\
x^2+7x,\quad 2x^2+7x,\quad 3x^2+7x,\\
x^2+9x,\quad 2x^2+9x,\quad 3x^2+11x.
\end{gathered}
\]

We derive upper bounds using combinatorial proofs and lower bounds using constructions of explicit colorings of $[N]$ found using the alternative computing inspired method of simulated annealing.

More generally, we also prove that the four color polynomial van der Waerden numbers for the entire family of nonzero polynomials \(ex^2+fx\) grow only linearly in the magnitude of their coefficients
\[
W(ex^2+fx;4)=\Theta(|e|+|f|),
\]
\end{abstract}

\section{Introduction}

\subsection{Van der Waerden numbers}

\begin{theorem}[van der Waerden~\cite{VanDerWaerden}]
For any integers $r,\ell\ge2$, there is an integer $N$ such that every
$r$ coloring of $[N]=\{1,\ldots,N\}$ contains
\[
a,\ a+d,\ \ldots,\ a+(\ell-1)d
\]
of one color, for some $d>0$.
\end{theorem}

Let $w(r,\ell)$ be the least such $N$. The upper bounds obtained from
van der Waerden's original proof were not primitive recursive. Shelah
\cite{ShelahClassical} improved on them by obtaining primitive recursive bounds, and Gowers
\cite{Gowers} later proved the explicit upper bound
\[
w(r,\ell)\le 2^{2^{r^{2^{2^{\ell+9}}}}}.
\]
Exponential lower bounds
are also known for fixed $r$ as $\ell$ grows \cite{GasarchHaeupler}.
The seven known nontrivial exact values are
\cite{KourilPaul,Kouril293,Monroe}:
\begingroup
\renewcommand{\arraystretch}{1.16}
\setlength{\arrayrulewidth}{0.4pt}
\setlength{\arraycolsep}{5pt}
\[
\begin{array}{|c|c|c|c|c|c|c|c|}\hline
(r,\ell) & (2,3) & (2,4) & (2,5) & (2,6) & (3,3) & (3,4) & (4,3) \\ \hline
w(r,\ell) & 9 & 35 & 178 & 1132 & 27 & 293 & 76 \\ \hline
\end{array}
\]
\endgroup
The known exact values are far smaller than the
general upper bounds. Our quadratic results give another example of
this gap and determine the dependence on the polynomial coefficients
up to constant factors.
\subsection{Polynomial van der Waerden numbers}

The polynomial van der Waerden theorem replaces equally spaced terms
by polynomial increments. Bergelson and Leibman proved the following theorem:

\begin{theorem}[Bergelson and Leibman~\cite{BergelsonLeibman}]
Let $P_1,\ldots,P_m\in\Z[x]$ satisfy $P_i(0)=0$. For any number of
colors $r$, every $r$ coloring of a sufficiently long interval $[N]$
contains
\[
a,\ a+P_1(d),\ \ldots,\ a+P_m(d)
\]
of one color for some positive integer $d$.
\end{theorem}

This original proof established the finiteness of these numbers but did not provide explicit upper bounds on their values. Walters \cite{Walters} subsequently gave a combinatorial proof with bounds that were not primitive recursive, while Shelah \cite{ShelahPolynomial} obtained primitive recursive bounds that remain extremely large. The classical van der Waerden theorem is itself a special case of the polynomial theorem. Taking the polynomials $x,2x,\ldots,(\ell-1)x$ gives an arithmetic progression of length $\ell$. Thus the large gap between exponential lower bounds and known upper bounds is also present in polynomial van der Waerden numbers which we define as follows:

\begin{definition}
For a nonzero $P(x)\in\Z[x]$ with $P(0)=0$ and $r\ge2$, the
\emph{$r$ color polynomial van der Waerden number} of $P$,
denoted $\PVDW{P(x)}{r}$, is the least positive integer $N$
such that every $r$ coloring of $[N]$ contains distinct
equally colored integers $u,v$ with $v-u=P(d)$ for some
$d\in\Z$.
\end{definition}

Equivalently, let
\[
\Dist{P}=\{\abs{P(d)}:d\in\Z,\ P(d)\ne0\}.
\]
For instance, $\Dist{x^2+x}=\{2,6,12,20,30,42,\ldots\}$.
Join $u,v\in[N]$ when $|u-v|\in\Dist{P}$. Then
$\PVDW{P(x)}{4}$ is the first $N$ for which this distance graph
requires at least five colors. Thus each exact value determines a
threshold for 4-colorability for an explicit family of finite graphs.

\subsection{Prior results}

Gasarch et al. \cite{GasarchEtAl} proved
\[
\begin{aligned}
\PVDW{x^2}{4}&\le84{,}149{,}474{,}894{,}213{,}522,
&\PVDW{x^2+x}{4}&\le12{,}828{,}393{,}907,\\
\PVDW{x^2+2x}{4}&\le66{,}287{,}711{,}296,
&\PVDW{x^2+3x}{4}&\le10{,}952{,}789.
\end{aligned}
\]
These bounds improve on the general proofs for the four quadratic polynomials shown above.
Canakci et al. \cite{CanakciEtAl} later used SAT solvers to prove
$\PVDW{x^2}{4}=58$.
\subsection{Our contribution}

We investigate simulated annealing as a method for discovering exact
polynomial van der Waerden numbers. For $x^2+x$, $2x^2+x$, and
$3x^2+x$, it found colorings of lengths $96$, $48$, and $32$.
Repeated attempts to extend these colorings were unsuccessful,
suggesting the exact values $97$, $49$, and $33$. We then proved
these conjectures using a single lemma involving $17$ integers.
This demonstrates how probabilistic search methods can guide exact proofs.
We also establish six further exact values by examining the
restrictions imposed by polynomial differences. For example, for
$2x^2+7x$, triangles in the distance graph force integers $12$
apart to have the same color. These constraints lead to a
contradiction on $[25]$.

Finally, we establish bounds that apply to every nonzero polynomial
$ex^2+fx$, regardless of the size of its coefficients:
\[
\frac49(|e|+|f|)+1
\le \PVDW{ex^2+fx}{4}
\le 64(|e|+|f|)+1.
\]
The main idea behind our upper bound is that we do not need a
separate computation for every polynomial. For almost all choices
of $e$ and $f$, a four coloring of a sufficiently long interval
would give a four coloring of the same graph on $715$ vertices.
We prove computationally that this graph cannot be four colored
and handle the exceptional cases separately. We also give methods for constructing four colorings that establish
lower bounds of the same asymptotic order. Consequently,
\[
\PVDW{ex^2+fx}{4}=\Theta(|e|+|f|).
\]
This shows that the true growth for this entire family is linear
in the coefficient sizes, despite the enormous general upper
bounds from the polynomial van der Waerden theorem.

\section{Exact Values}\label{sec:exact}

We use the following identity from Lemma~4 of Gasarch
et al.\ \cite{GasarchEtAl}.

\begin{proposition}\label{prop:multiplication}
For every positive integer $h$ and every $r\ge2$,
\[
\PVDW{hP(x)}{r}=h\bigl(\PVDW{P(x)}{r}-1\bigr)+1.
\]
\end{proposition}

The same elementary fact is used in the first three upper bounds.
\begin{lemma}\label{lem:17}
Every four coloring of
\[
\{0,1,\ldots,16\}
\]
contains two equally colored integers whose difference belongs to
\[
\mathcal F=\{1,2,5,7,12,15\}.
\]
\end{lemma}

\begin{proof}
Suppose, for contradiction, that there is a four coloring of
\[
\{0,1,\ldots,16\}
\]
with no two equally colored integers at a difference in $\mathcal F$.

There are $17$ integers and only four colors. If every color were used at
most four times, then at most $4\cdot4=16$ integers would be colored.
Hence some color is used at least five times. Choose five integers of that
color and write them in increasing order as
\[
a_0<a_1<a_2<a_3<a_4.
\]

For $1\le i\le4$, let
$g_i=a_i-a_{i-1}$.

By our assumption, no two equally colored integers may differ by an
element of $\mathcal F$. Since $a_{i-1}$ and $a_i$ have the same color, we must have
\[
g_i\notin\{1,2,5,7\}.
\]
In particular, $g_i\ge3$.

Also,
\[
g_1+g_2+g_3+g_4=a_4-a_0\le16.
\]
If some $g_i\ge8$, then the other three gaps are each at least $3$, so
\[
g_1+g_2+g_3+g_4\ge8+3+3+3=17,
\]
a contradiction. Therefore every $g_i\le7$.

Hence each gap is an integer between $3$ and $7$, but it cannot equal
$5$ or $7$. Therefore
\[
g_i\in\{3,4,6\}
\qquad(1\le i\le4).
\]
We now consider how many of the four gaps are equal to $6$. Suppose first that no gap is $6$. Then every gap is either $3$ or $4$.
A $3$ and a $4$ cannot occur next to one another, since two chosen
positions would then differ by
$3+4=7$,
a forbidden distance. Hence all four gaps are equal. If they are all $3$,
then
$a_4-a_0=3+3+3+3=12$,
while if they are all $4$, then
$a_3-a_0=4+4+4=12$.
Both cases give two equally colored positions at the forbidden difference
$12$.

Next, there cannot be two gaps equal to $6$, since the remaining two gaps
are at least $3$, giving
$g_1+g_2+g_3+g_4\ge6+6+3+3=18>16$.

It remains to consider exactly one gap equal to $6$. If the other three
gaps are all $3$, then
$a_4-a_0=6+3+3+3=15$,
a forbidden distance. Otherwise one of the remaining gaps is $4$. There
cannot be two gaps equal to $4$, since then
$6+4+4+3=17>16$.
Thus the four gaps are $6,4,3,3$ in some order. Since a gap $4$ cannot
be adjacent to a gap $3$, the $4$ must be next to the $6$ and at one end
of the sequence. Up to reversing the order, the gaps are therefore
\[
4,6,3,3.
\]
The final three gaps sum to
$6+3+3=12$,
so again two equally colored chosen positions differ by the forbidden
distance $12$.

All possibilities lead to a contradiction. This proves the lemma.
\end{proof}

Lemma~\ref{lem:17} gives a sufficient condition for a short exact proof.
If $h$ is a positive integer with $h\mathcal F\subseteq\Dist{P}$,
the colors of $1,1+h,\ldots,1+16h$ would contradict the lemma.
Thus $\PVDW{P(x)}{4}\le16h+1$.
A valid four coloring of $[16h]$ then proves equality.
The next three theorems use $h=6,3,2$, respectively.

\begin{theorem}\label{thm:x2x}
\[
\boxed{\PVDW{x^2+x}{4}=97.}
\]
\end{theorem}

\begin{proof}
\textbf{Upper bound.}
Suppose that $[97]$ has a valid four coloring for $P(x)=x^2+x$. Consider the $17$ integers
\[
1,7,13,\ldots,97.
\]
These are exactly the numbers $1+6j$ for $0\le j\le16$. Give the index $j$ the color of $1+6j$. By Lemma~\ref{lem:17}, there are indices $i<j$ with the same color and
\[
j-i\in\mathcal F.
\]
Therefore the corresponding integers in $[97]$ have the same color and differ by
\[
6(j-i)\in6\mathcal F=\{6,12,30,42,72,90\}.
\]
Each of these is a value of $d^2+d$:
\[
\begin{gathered}
6=2^2+2,\quad12=3^2+3,\quad30=5^2+5,\\
42=6^2+6,\quad72=8^2+8,\quad90=9^2+9.
\end{gathered}
\]
Thus the assumed coloring contains an equally colored pair at a forbidden polynomial distance, a contradiction. Hence
\[
\PVDW{x^2+x}{4}\le97.
\]

\textbf{Lower bound.}
The nonzero distances in $\Dist{x^2+x}$ that are smaller than $97$ are
\[
2,6,12,20,30,42,56,72,90.
\]
The following coloring of $[96]$ was found using simulated annealing. For each of the nine forbidden distances listed above, comparing every pair of entries at that separation shows that the two entries have different colors.
\witnessbox{\Rcol\Rcol\Gcol\Gcol\Rcol\Rcol\Ycol\Ycol\Rcol\Rcol\Ycol\Ycol\Gcol\Gcol\Bcol\Bcol\Gcol\Gcol\Rcol\Rcol\Gcol\Gcol\Rcol\Rcol\Bcol\Bcol\Rcol\Rcol\Bcol\Bcol\Gcol\Gcol\Ycol\Ycol\Gcol\Gcol\Rcol\Rcol\Bcol\Bcol\Rcol\Rcol\Bcol\Bcol\Rcol\Rcol\Bcol\Bcol\\[-1pt]
\Gcol\Gcol\Ycol\Ycol\Gcol\Gcol\Ycol\Ycol\Gcol\Gcol\Ycol\Ycol\Bcol\Bcol\Rcol\Rcol\Bcol\Bcol\Gcol\Gcol\Ycol\Ycol\Gcol\Gcol\Ycol\Ycol\Bcol\Bcol\Ycol\Ycol\Bcol\Bcol\Gcol\Gcol\Bcol\Bcol\Rcol\Rcol\Ycol\Ycol\Rcol\Rcol\Ycol\Ycol\Bcol\Bcol\Ycol\Ycol}
Hence a valid four coloring exists on $[96]$, so
\[
\PVDW{x^2+x}{4}\ge97.\qedhere
\]
\end{proof}

\begin{theorem}\label{thm:2x2x}
\[\boxed{\PVDW{2x^2+x}{4}=49.}\]
\end{theorem}
\begin{proof}
\textbf{Upper bound.} Apply Lemma~\ref{lem:17} to the colors of
$1+3j$ for $0\le j\le16$, as in the preceding proof.
An equally colored pair has difference in
\[3\mathcal F=\{3,6,15,21,36,45\}.\]
For $P(x)=2x^2+x$, these are the values at
$d=1,-2,-3,3,4,-5$, respectively. They are all forbidden, so
$\PVDW{2x^2+x}{4}\le49$.

\textbf{Lower bound.} The values
\[
P(m)=\frac{(2m)(2m+1)}2,\qquad
P(-m)=\frac{(2m-1)(2m)}2\qquad(m\ge1)
\]
give every positive triangular number. The forbidden distances below
$49$ are $1,3,6,10,15,21,28,36,45$.
The following coloring of $[48]$ was found using
simulated annealing. Direct comparison at each listed distance
shows that no forbidden pair has the same color.
\witnessbox{\Rcol\Gcol\Rcol\Ycol\Rcol\Ycol\Gcol\Bcol\Gcol\Rcol\Gcol\Rcol\Bcol\Rcol\Bcol\Gcol\Ycol\Gcol\Rcol\Bcol\Rcol\Bcol\Rcol\Bcol\Gcol\Ycol\Gcol\Ycol\Gcol\Ycol\Bcol\Rcol\Bcol\Gcol\Ycol\Gcol\Ycol\Bcol\Ycol\Bcol\Gcol\Bcol\Rcol\Ycol\Rcol\Ycol\Bcol\Ycol}
Thus $\PVDW{2x^2+x}{4}\ge49$.
\end{proof}

\begin{theorem}\label{thm:3x2x}
\[\boxed{\PVDW{3x^2+x}{4}=33.}\]
\end{theorem}
\begin{proof}
\textbf{Upper bound.} This time apply Lemma~\ref{lem:17} to the colors
of $1+2j$ for $0\le j\le16$. The resulting difference belongs to
\[2\mathcal F=\{2,4,10,14,24,30\}.\]
For $P(x)=3x^2+x$, these values occur at
$d=-1,1,-2,2,-3,3$. Hence \mbox{$\PVDW{3x^2+x}{4}\le33$}.

\textbf{Lower bound.} The same six values are all the forbidden distances
below $33$, since $|3d^2+d|\ge3d^2-|d|\ge44$ for $|d|\ge4$.
The following coloring of $[32]$ was found using
simulated annealing. Entries at each of the six listed distances
have different colors.
\witnessbox{\Gcol\Gcol\Bcol\Bcol\Rcol\Rcol\Gcol\Gcol\Ycol\Ycol\Rcol\Rcol\Gcol\Gcol\Ycol\Ycol\Rcol\Rcol\Bcol\Bcol\Ycol\Ycol\Rcol\Rcol\Bcol\Bcol\Ycol\Ycol\Gcol\Gcol\Bcol\Bcol}
Therefore $\PVDW{3x^2+x}{4}\ge33$.
\end{proof}

\subsection{Odd cycles}

Join two integers when their distance is forbidden. The next lemma is
used in both upper bound proofs in this subsection.
\begin{lemma}\label{lem:oddcycle}
In a proper four coloring, two vertices adjacent to every vertex of the
same odd cycle have the same color.
\end{lemma}
\begin{proof}
The cycle needs at least three colors. Each additional vertex excludes
its color from the cycle, so both must use the remaining fourth color.
\end{proof}

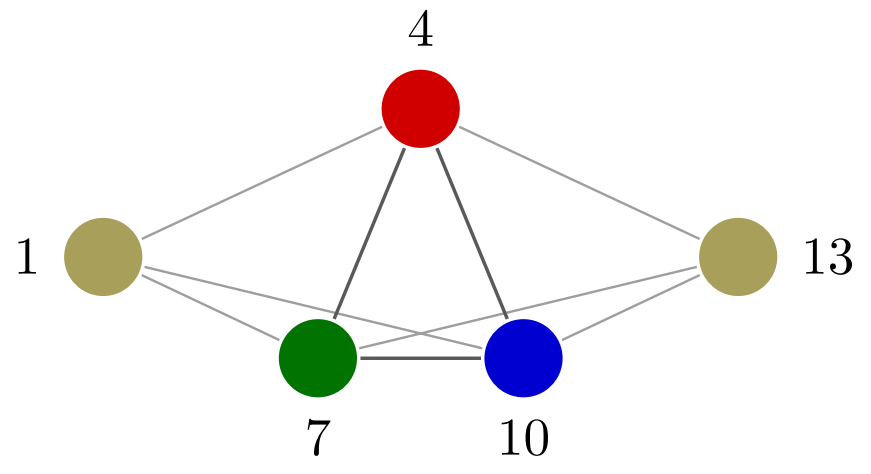
\begin{figure}[htbp]
\centering
\begin{tikzpicture}[font=\small]
\coordinate (a) at (0,1.0);
\coordinate (b) at (-.68,-.65);
\coordinate (c) at (.68,-.65);
\coordinate (p) at (-2.1,.02);
\coordinate (q) at (2.1,.02);
\draw[witness edge] (p)--(a) (p)--(b) (p)--(c)
                     (q)--(a) (q)--(b) (q)--(c);
\draw[black!65,line width=.65pt] (a)--(b)--(c)--cycle;
\node[witness point,fill=witnessred,label=above:4] at (a) {};
\node[witness point,fill=witnessgreen,label=below:7] at (b) {};
\node[witness point,fill=witnessblue,label=below:10] at (c) {};
\node[witness point,fill=witnessyellow,label=left:1] at (p) {};
\node[witness point,fill=witnessyellow,label=right:13] at (q) {};
\end{tikzpicture}
\caption{A common triangle for $2x^2+7x$. Every drawn edge has forbidden
length $3$, $6$, or $9$. The vertices $4,7,10$ use three colors, so
$1$ and $13$ have the same fourth color.}
\label{fig:forcing}
\end{figure}

\begin{theorem}\label{thm:seven}
\[\boxed{\PVDW{2x^2+7x}{4}=25.}\]
\end{theorem}
\begin{proof}
\textbf{Upper bound.} The forbidden distances below $25$ are
$3,4,5,6,9,15,22$. Suppose $c$ is a valid four coloring of $[25]$.
For $1\le i\le13$, the vertices $i+3,i+6,i+9$ form a triangle.
Both $i$ and $i+12$ are adjacent to all three. Lemma~\ref{lem:oddcycle} gives
\[c(i)=c(i+12)\qquad(1\le i\le13).\]
Arrange the positions $1,\ldots,12$ around a circle. Two positions separated
by three, four, five, or six steps around this circle have different colors.
Indeed, suitable representatives in $[24]$ differ by one of $3,4,5,6$.
Thus each color appears at most three times, and three positions of one color
must be consecutive around the circle.

The color of position $1$ cannot appear three times. The possible triples
would be $(11,12,1)$, $(12,1,2)$, and $(1,2,3)$. They are ruled out,
respectively, by the forbidden pairs $(1,23)$, $(2,24)$, and $(3,25)$,
each at distance $22$, together with the equalities already proved.
Hence the four colors cover at most $2+3+3+3=11$ of the twelve positions,
a contradiction.

\Needspace{6\baselineskip}
\textbf{Lower bound.} The following coloring of $[24]$ was recovered by
simulated annealing. Direct comparison at the seven listed distances gives
no monochromatic forbidden pair.
\witnessbox{\Ycol\Ycol\Ycol\Rcol\Rcol\Rcol\Gcol\Gcol\Gcol\Bcol\Bcol\Bcol\Ycol\Ycol\Ycol\Rcol\Rcol\Rcol\Gcol\Gcol\Gcol\Bcol\Bcol\Bcol}
Thus $\PVDW{2x^2+7x}{4}\ge25$.
\end{proof}

\begin{theorem}\label{thm:nine}
\[\boxed{\PVDW{2x^2+9x}{4}=29.}\]
\end{theorem}
\begin{proof}
\textbf{Upper bound.} The forbidden distances below $29$ are
$4,5,7,9,10,11,18,26$. Suppose $c$ is a valid four coloring of $[29]$.
For $5\le i\le11$, the vertices
\[i-4,\quad i+5,\quad i+9,\quad i+18,\quad i+7\]
form a five cycle in the displayed order. Its edge lengths are
$9,4,9,11,11$. Both $i$ and $i+14$ are adjacent to every vertex of this
cycle, so Lemma~\ref{lem:oddcycle} gives
\[c(i)=c(i+14)\qquad(5\le i\le11).\]

Now let $1\le i\le4$ and consider the colors at
\[i+4,\quad i+7,\quad i+10,\quad i+5,\quad i+9.\]
The first two consecutive pairs have different colors because
$c(i+4)=c(i+18)$ and $c(i+7)=c(i+21)$, while the pairs
$(i+7,i+18)$ and $(i+10,i+21)$ have forbidden distance $11$.
The other three consecutive pairs have forbidden distances $5,4,5$.
These inequalities form a five cycle, so at least three colors occur. Both $i$ and $i+14$
are adjacent to all five displayed vertices, so their colors agree.
Reflection by $v\mapsto30-v$ gives $12\le i\le15$.
Thus
\[c(i)=c(i+14)\qquad(1\le i\le15).\]

Finally, consider the colors of $1,3,6,10,13$. The pairs
$(1,6)$, $(1,10)$, $(3,10)$, $(3,13)$, $(6,10)$, and $(6,13)$ have
different colors directly from their forbidden distances. The remaining four
inequalities follow from the established equalities and these edges:
\begin{center}
\begingroup
\small
\renewcommand{\arraystretch}{1.16}
\setlength{\arrayrulewidth}{0.4pt}
\setlength{\tabcolsep}{6pt}
\begin{tabular}{|c|c|c|}\hline
\textbf{Colors compared} & \textbf{Forbidden pair} & \textbf{Distance}\\\hline
$c(1),c(3)$ & $29,3$ & $26$\\\hline
$c(1),c(13)$ & $1,27$ & $26$\\\hline
$c(3),c(6)$ & $17,6$ & $11$\\\hline
$c(10),c(13)$ & $24,13$ & $11$\\\hline
\end{tabular}
\endgroup
\end{center}
Five distinct colors are required, a contradiction.

\Needspace{7\baselineskip}
\textbf{Lower bound.} The following annealed coloring of $[28]$ has no
monochromatic pair at any of the eight listed distances.
\witnessbox{\Ycol\Ycol\Ycol\Rcol\Rcol\Rcol\Rcol\Bcol\Bcol\Bcol\Gcol\Gcol\Gcol\Gcol\Ycol\Ycol\Ycol\Ycol\Rcol\Rcol\Rcol\Bcol\Bcol\Bcol\Bcol\Gcol\Gcol\Gcol}
Thus $\PVDW{2x^2+9x}{4}\ge29$.
\end{proof}

\Needspace{10\baselineskip}
\begin{proposition}\label{prop:oddtransfer}
For every odd integer $f$ and every $r\ge2$,
\[\PVDW{x^2+fx}{r}=2\PVDW{2x^2+fx}{r}-1.\]
\end{proposition}
\begin{proof}
Set $P(x)=x^2+fx$ and $Q(x)=2x^2+fx$.
For even $d=2t$, $P(d)=2Q(t)$. For odd $d$, the integer
$t=-(d+f)/2$ also satisfies $P(d)=2Q(t)$.
Conversely, $2Q(t)=P(2t)$. Thus $\Dist{P}=2\Dist{Q}$, and
Proposition~\ref{prop:multiplication} gives the result.
\end{proof}
In particular, Theorems~\ref{thm:seven} and \ref{thm:nine} give
\[
\boxed{\PVDW{x^2+7x}{4}=49,\qquad
\PVDW{x^2+9x}{4}=57.}
\]

\subsection{Block colorings}

\begin{figure}[htbp]
\centering
\begin{tikzpicture}[
  vertex/.style={
    circle,
    draw=white,
    line width=.7pt,
    minimum size=6.5mm,
    inner sep=0pt,
    font=\small\bfseries
  },
  edge/.style={draw=black!35,line width=.45pt}
]

\foreach \shift/\offset in {-1.65/0,1.65/1}{
  \begin{scope}[xshift=\shift cm]
    \foreach \k in {0,...,7}{
      \pgfmathtruncatemacro{\n}{2*\k+1+\offset}
      \coordinate (v\n) at ({90-45*\k}:1.28);
    }
  \end{scope}
}

\foreach \i in {1,...,16}{
  \foreach \d in {4,6,8,10,14}{
    \pgfmathtruncatemacro{\j}{\i+\d}
    \ifnum\j<17
      \draw[edge] (v\i)--(v\j);
    \fi
  }
}

\foreach \offset in {0,1}{
  \foreach \k/\colour/\ink in {
    0/red!82!black/white,
    1/red!82!black/white,
    2/green!45!black/white,
    3/green!45!black/white,
    4/blue!82!black/white,
    5/blue!82!black/white,
    6/yellow!58!black/black,
    7/yellow!58!black/black%
  }{
    \pgfmathtruncatemacro{\n}{2*\k+1+\offset}
    \node[vertex,fill=\colour,text=\ink] at (v\n) {\n};
  }
}

\end{tikzpicture}
\caption{A valid coloring of $[16]$ for $3x^2+11x$.
The $38$ edges have lengths $4,6,8,10,14$, and every edge has
differently colored endpoints. All forbidden distances are even,
so the odd and even vertices form separate components.}
\label{fig:witness17}
\end{figure}
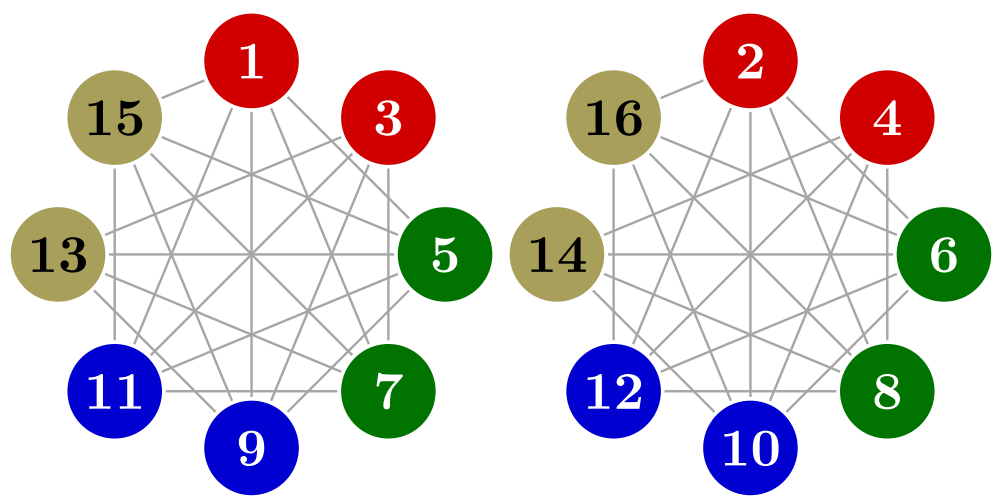

\begin{theorem}\label{prop:smallblockexact}
\[\boxed{\PVDW{3x^2+7x}{4}=41,\qquad
\PVDW{3x^2+11x}{4}=17.}\]
\end{theorem}
\begin{proof}
\textbf{Upper bound.} Restrict a coloring of $[41]$ to its $21$ odd
integers, written as $1+2j$ for $0\le j\le20$.
The distances $1,2,3,20$ between the corresponding values of $j$ are
forbidden, since $2,4,6,40\in\Dist{3x^2+7x}$.
Six values of $j$ of one color would have five consecutive gaps of at
least $4$. Their total is at most $20$, so all five gaps would be $4$.
The first and last would differ by the forbidden distance $20$.
Thus each color covers at most five of the $21$ positions, a contradiction.

\textbf{Lower bound.} Repeat the four blocks \Rcol\Rcol,
\Gcol\Gcol, \Bcol\Bcol, \Ycol\Ycol\ five times:
\witnessbox{\Rcol\Rcol\Gcol\Gcol\Bcol\Bcol\Ycol\Ycol\Rcol\Rcol\Gcol\Gcol\Bcol\Bcol\Ycol\Ycol\Rcol\Rcol\Gcol\Gcol\Bcol\Bcol\Ycol\Ycol\Rcol\Rcol\Gcol\Gcol\Bcol\Bcol\Ycol\Ycol\Rcol\Rcol\Gcol\Gcol\Bcol\Bcol\Ycol\Ycol}
The forbidden distances below $41$ are $2,4,6,10,20,26,40$.
An equally colored pair in this word has distance $1$ or a distance in
$[7,9]$, $[15,17]$, $[23,25]$, or $[31,33]$.
None is forbidden, so the coloring of $[40]$ is valid.

\textbf{Upper bound.} Restrict $[17]$ to its nine odd integers.
After dividing their differences by two, the forbidden distances are
$2,3,4,5,7$. Three positions of one color would have two consecutive gaps
from $\{1,6,8\}$ with sum at most $8$.
The only possibilities are $(1,1)$, $(1,6)$, and $(6,1)$.
Their sums $2$ and $7$ are forbidden. Each color therefore covers at most
two of the nine positions, again a contradiction.

\textbf{Lower bound.} The smallest positive distance in $\Dist{3x^2+11x}$ is $4$.
Four constant blocks of length $4$ give a valid coloring of $[16]$:
\witnessbox{\Rcol\Rcol\Rcol\Rcol\Gcol\Gcol\Gcol\Gcol\Bcol\Bcol\Bcol\Bcol\Ycol\Ycol\Ycol\Ycol}
Every equally colored pair has distance at most $3$.
Figure~\ref{fig:witness17} shows all the forbidden pairs.
\end{proof}

\FloatBarrier
\section{Asymptotic Lower Bounds}\label{sec:lower}

The exact values above depend on the forbidden distances of individual
polynomials. We now construct valid colorings that give lower bounds
without solving each coloring problem separately.

\begin{proposition}\label{prop:shortlower}
Let $s_1<s_2$ be the two smallest members of $\Dist{P}$. Then
\[\boxed{\PVDW{P(x)}{4}\ge\max\{4s_1,2s_2\}+1.}\]
\end{proposition}
\begin{proof}
Color four consecutive blocks of length $s_1$ with four different colors.
An equally colored pair lies in one block and has distance less than $s_1$.
This gives a valid coloring of $[4s_1]$.

For the other bound, split $[2s_2]$ into two intervals of length $s_2$.
Within either interval, the only possible forbidden distance is $s_1$.
These edges form disjoint paths. Color each path alternately with two
colors. Use a different pair of colors on the second interval.
\end{proof}

\begin{proposition}\label{prop:blocklower}
Suppose $h,J\ge1$ are integers, no forbidden distance is smaller than $h$,
and, for $1\le j<J$, no forbidden distance lies in
\[[(4j-1)h+1,(4j+1)h-1].\]
Then
\[\PVDW{P(x)}{4}\ge4hJ+1.\]
\end{proposition}
\begin{proof}
Repeat $J$ times the word consisting of $h$ copies of each of
\Rcol, \Gcol, \Bcol, \Ycol, in that order.
An equally colored pair in the same block has distance at most $h-1$.
If its blocks are $4j$ blocks apart, its distance lies in the displayed
interval. Both possibilities are excluded by the hypotheses.
\end{proof}
The choices $(h,J)=(3,2),(2,5),(4,1)$ give the exact lower bounds
for $2x^2+7x$, $3x^2+7x$, and $3x^2+11x$, respectively.

\subsection{Lower bounds in terms of the coefficients}

The first constructions depend on the smallest forbidden distances.
For quadratic polynomials we can also give lower bounds directly in terms of the
coefficients, using the following consequence of a theorem of Barajas
and Serra \cite[Theorem~5]{BarajasSerra}.
\begin{lemma}\label{lem:fourdistances}
Let $s_1<s_2<s_3<s_4$ be positive integers with $s_3>s_1+s_2$.
The integers have a four coloring avoiding all four distances.
\end{lemma}
\begin{proof}
After division by their greatest common divisor, every set of four
distances requiring five colors has its third smallest distance equal
to the sum of the first two \cite[Theorem~5]{BarajasSerra}.
The strict inequality excludes this case. A coloring for the normalized
distances extends along each residue class of the common divisor.
\end{proof}

\begin{theorem}\label{thm:regionlower}
For $E,F>0$, the following lower bounds hold:
\[
\begin{array}{ll}
\PVDW{Ex^2+Fx}{4}\ge9E-3F+1,&0<F<E,\\[3pt]
\PVDW{Ex^2+Fx}{4}\ge4E+2F+1,&E<F<2E,\\[3pt]
\PVDW{Ex^2+Fx}{4}\ge\left\lceil\dfrac{4(F-E)}3\right\rceil+1,&F\ge E.
\end{array}
\]
\end{theorem}
\begin{proof}
If $0<F<E$, the first five distinct forbidden distances are
\[E-F,\ E+F,\ 4E-2F,\ 4E+2F,\ 9E-3F.\]
The third exceeds the sum of the first two. Lemma~\ref{lem:fourdistances}
therefore gives a four coloring avoiding the first four. Restrict it to
$[9E-3F]$, where no other forbidden distance can occur.

If $E<F<2E$, write $\delta=F-E$, so $0<\delta<E$. Below $4E+2F=6E+2\delta$, the only
possible forbidden distances are
\[\delta,\quad2E-2\delta,\quad2E+\delta,\quad6E-3\delta.\]
Indeed, the next positive input value is $6E+2\delta$, and the next
negative input value is $12E-4\delta>6E+2\delta$.
The first two displayed quantities are the smallest. If they coincide,
there are at most three distinct distances, and greedy coloring uses four
colors. Otherwise, the third exceeds their sum. Apply
Lemma~\ref{lem:fourdistances} and restrict to $[4E+2F]$.

For the final bound, write $F=E(n+\rho)$ with $n\ge1$ an integer and
$0\le\rho<1$. All positive forbidden distances except possibly
\[En\rho\quad\text{and}\quad E(n+1)(1-\rho)\]
are at least $F-E$. To check the negative inputs $-m$ with $1\le m<n$,
the function $m(F-Em)$ is concave. Its endpoint values are $F-E$ and
$E(n-1)(1+\rho)$, both at least $F-E$. When $n=1$ this range is empty.
Inputs with $m\ge n+2$ give larger values, while positive inputs give at
least $E+F$.

The only possible positive forbidden distances smaller than $F-E$ therefore arise from $d=-n$ and $d=-(n+1)$. Their absolute values sum to
\[
En\rho+E(n+1)(1-\rho)=E(n+1-\rho)\ge F-E.
\]
Let $s_1<s_2$ be the two smallest positive forbidden distances.
If $s_1\ge(F-E)/3$, then $4s_1\ge4(F-E)/3$.
Otherwise, $s_1<(F-E)/3$ must be one of the two exceptional
distances above. Since their sum is at least $F-E$ and
all other forbidden distances are at least $F-E$,
we have $s_2\ge(F-E)-s_1>2(F-E)/3$.

Thus, in either case,
\[
\max\{4s_1,2s_2\}\ge\frac{4(F-E)}3.
\]
Proposition~\ref{prop:shortlower} now gives
\[
\PVDW{Ex^2+Fx}{4}
\ge\left\lceil\frac{4(F-E)}3\right\rceil+1.
\]
\end{proof}
\begin{corollary}
For all integers $e,f$, not both zero,
\[
W(ex^2+fx;4)\ge\frac49(|e|+|f|)+1.
\]
\end{corollary}

\begin{proof}
By symmetry, assume $E=|e|$ and $F=|f|$.

If $0<F<E$, Theorem~\ref{thm:regionlower} gives
\[
W(Ex^2+Fx;4)\ge9E-3F+1
>\frac49(E+F)+1.
\]

If $E<F<2E$, the same theorem gives
\[
W(Ex^2+Fx;4)\ge4E+2F+1
>\frac49(E+F)+1.
\]

If $F\ge2E$, then
\[
\frac43(F-E)\ge\frac49(E+F),
\]
so the third bound of Theorem~\ref{thm:regionlower}
implies the result.

Finally, if $E=0$, $F=0$, or $E=F$,
the smallest positive forbidden distance is,
respectively, $F$, $E$, or $2E$.
Proposition~\ref{prop:shortlower} gives the
claimed bound in each case.
\end{proof}

\section{Asymptotic Upper Bound}\label{sec:upper}

For $P(x)=Ex^2+Fx$ with $E,F\ge0$, label a lattice point $(u,v)$ by
$\phi(u,v)=1+Fu+Ev$. Then
\[
\begin{aligned}
\phi(u+d,v+d^2)-\phi(u,v)&=Ed^2+Fd=P(d),\\
\phi(u+d,v-d^2)-\phi(u,v)&=Fd-Ed^2=-P(-d).
\end{aligned}
\]
Thus pairs separated by $(d,\pm d^2)$ give polynomial differences.
The following finite graph is used for every choice of the coefficients.

\begin{lemma}\label{lem:grid}
Let $Q=\{0,1,\ldots,10\}\times\{0,1,\ldots,64\}$.
For $1\le d\le8$, join $(u,v)$ to $(u+d,v+d^2)$ and to
$(u+d,v-d^2)$ whenever both points lie in $Q$.
Every four coloring of $Q$ contains two equally colored joined points.
\end{lemma}
\begin{proof}
The graph has $715$ vertices and $4864$ edges.
The exact search finds no proper four coloring. The search procedure is
described in Section~\ref{sec:methods}.
\end{proof}

An edge of displacement $(d,-d^2)$ has label difference $d(F-Ed)$.
If $F=dE$, its two endpoints receive the same integer label, so this edge
cannot give a nonzero polynomial difference. Since $1\le d\le8$, the
only excluded cases are $F=E,2E,\ldots,8E$.

\begin{lemma}\label{lem:exceptional}
For every $1\le j\le8$, $\PVDW{x^2+jx}{4}\le97$.
\end{lemma}
\begin{proof}
The case $j=1$ is Theorem~\ref{thm:x2x}.
For $j=2,\ldots,8$, the exact search is applied to $[97]$,
joining two integers when their distance is in $\Dist{x^2+jx}$.
In each of the seven cases the search finds no proper four coloring.
\end{proof}

\begin{theorem}\label{thm:main}
For all integers $e,f$, not both zero,
\[\boxed{\PVDW{ex^2+fx}{4}\le64\bigl(|e|+|f|\bigr)+1.}\]
\end{theorem}
\begin{proof}
Set $E=|e|$ and $F=|f|$. Replacing $P(x)$ by $-P(x)$ or $P(-x)$
does not change its forbidden distances, so assume $P(x)=Ex^2+Fx$.
If $E=0$, the identity $\PVDW{Fx}{4}=4F+1$ from
Gasarch et al.\ \cite{GasarchEtAl} gives the result. Hence assume $E>0$.

First suppose $F\notin\{E,2E,\ldots,8E\}$.
A valid coloring of $[64E+10F+1]$ would color every $(u,v)\in Q$ by
the color of $1+Fu+Ev$. Lemma~\ref{lem:grid} gives equally colored
points separated by $(d,d^2)$ or $(d,-d^2)$ for some $1\le d\le8$.
Their integer labels differ by $P(d)$ or $-P(-d)$, respectively.
The first is positive and the second is nonzero because $F\ne dE$.
Either gives a forbidden equally colored pair, a contradiction. Thus
\[
\PVDW{Ex^2+Fx}{4}\le64E+10F+1\le64(E+F)+1.
\]

If $F=jE$ for $1\le j\le8$, Lemma~\ref{lem:exceptional} and
Proposition~\ref{prop:multiplication} give
\[
\PVDW{E(x^2+jx)}{4}\le96E+1\le64(E+F)+1,
\]
since $j\ge1$. This completes the proof.
\end{proof}

The known value
$\PVDW{x^2}{4}=58$ \cite{CanakciEtAl} and
Proposition~\ref{prop:multiplication} give $\PVDW{Ex^2}{4}=57E+1$.
Consequently, if $C_4$ is the smallest constant such that
$\PVDW{ex^2+fx}{4}\le C_4(|e|+|f|)+1$ for all nonzero pairs, then
\[57\le C_4\le64.\]

Together with the coefficient lower bounds in Section~\ref{sec:lower},
this establishes
$\PVDW{ex^2+fx}{4}=\Theta(|e|+|f|)$
for all nonzero integer coefficient pairs $(e,f)$.

\subsection{Comparison of Bounds to Exact Values}
Our upper bound $64(|e|+|f|)+1$ applies to every
nonzero polynomial $ex^2+fx$.
Table~\ref{tab:comparison} compares this bound with
our nine exact values.

\begin{table}[htbp]
\centering
\begingroup
\small
\renewcommand{\arraystretch}{1.16}
\setlength{\arrayrulewidth}{0.4pt}
\setlength{\tabcolsep}{8pt}
\begin{tabular}{|l|r|r|}\hline
\textbf{Polynomial} & \textbf{Exact} & \textbf{UB}\\\hline
$x^2+x$ & 97 & 129\\\hline
$2x^2+x$ & 49 & 193\\\hline
$3x^2+x$ & 33 & 257\\\hline
$2x^2+7x$ & 25 & 577\\\hline
$x^2+7x$ & 49 & 513\\\hline
$2x^2+9x$ & 29 & 705\\\hline
$x^2+9x$ & 57 & 641\\\hline
$3x^2+7x$ & 41 & 641\\\hline
$3x^2+11x$ & 17 & 897\\\hline
\end{tabular}
\endgroup
\caption{Exact values and uniform upper bounds for our nine quadratic polynomials}
\label{tab:comparison}
\end{table}

For these nine polynomials, the upper
bounds are within a factor of $53$. Despite applying to every
nonzero polynomial $ex^2+fx$, these bounds provide reasonably
close estimates in several cases. For example, our upper bound
for $x^2+x$ is $129$, only $32$ greater than the exact value $97$.

The improvement over earlier explicit upper bounds is much
greater. Gasarch et al.\ obtained upper bounds of
$12{,}828{,}393{,}907$ for $x^2+x$ and approximately
$6.57\times10^{18}$ for $3x^2+x$. Our bound reduces
these to $129$ and $257$, respectively. Thus, our bounds not
only establish the correct linear dependence on the coefficients
but also give numerically useful estimates that are substantially
closer to the exact values than earlier general estimates.

\FloatBarrier

\section{Computational Methods}\label{sec:methods}

\subsection{Simulated annealing}

For a polynomial $P(x)$ and interval $[N]$, let $E_P(N)$ be the set of forbidden pairs. For a four coloring
$c:[N]\to\{1,2,3,4\}$, define
\[
E(c)=
\#\{\{u,v\}\in E_P(N):c(u)=c(v)\}.
\]
Thus $E(c)$ is the number of monochromatic forbidden pairs, and $E(c)=0$
exactly when $c$ is a valid four coloring.

The energy change from recoloring one vertex depends only on its
neighbors. We use the Metropolis simulated annealing rule
\cite{KirkpatrickEtAl}: a recoloring that does not increase $E$ is
accepted, while one increasing it by $\Delta>0$ is accepted with
probability $e^{-\Delta/T}$. We decrease $T$ geometrically during
each run.
\begin{algorithm}[H]
\caption{Simulated annealing for the lower bounds}
\label{alg:anneal}
\begin{algorithmic}[1]
\Require Forbidden pairs $E_P(N)$, restarts $R$, proposals $L$, temperatures $T_0,T_{\min}$
\State Set $\alpha\gets(T_{\min}/T_0)^{1/L}$
\For{$j=1,\ldots,R$}
    \State Choose a random four coloring $c$ and set $T\gets T_0$
    \If{$E(c)=0$}
        \State \Return $c$
    \EndIf
    \For{$s=1,\ldots,L$}
        \State Choose a vertex $v$ and a color $q\ne c(v)$ uniformly at random
        \State Let $\Delta$ be the energy change from recoloring $v$ by $q$, and draw $U$ uniformly from $[0,1]$
        \If{$\Delta\le0$ or $U<e^{-\Delta/T}$}
            \State Set $c(v)\gets q$
        \EndIf
        \If{$E(c)=0$}
            \State \Return $c$
        \EndIf
        \State Set $T\gets\alpha T$
    \EndFor
\EndFor
\State \Return \textsc{failure}
\end{algorithmic}
\end{algorithm}

Simulated annealing has previously been applied to graph coloring
\cite{JohnsonEtAlAnnealing} and classical Ramsey problems
\cite{ExooTatarevic}. We adapt it to polynomial distance
graphs, where the objective is to find a four coloring with
$E(c)=0$. The difficulty is that a coloring can have positive
energy even when every single vertex recoloring has
$\Delta\ge0$. The Metropolis rule \cite{KirkpatrickEtAl}
allows the search to move away from such local minima, while
decreasing the temperature gradually favors colorings with fewer
forbidden monochromatic pairs.
We use independent random restarts and search for valid colorings
on successively longer intervals. Kissat \cite{Kissat} and other
deterministic search methods struggled on some of these instances,
motivating our use of this probabilistic approach. Whenever the
search finds a coloring with $E(c)=0$, we verify it by checking
every pair in $E_P(N)$. This establishes the rigorous lower bound
$\PVDW{P(x)}{4}\ge N+1$.
\FloatBarrier
\subsection{Exact search}

The same exact search is used for both Lemma~\ref{lem:grid} and
Lemma~\ref{lem:exceptional}.

The search uses the saturation degree rule of Br\'elaz \cite{Brelaz}.
At each step it chooses an uncolored vertex having the largest number of
different colors among its colored neighbors. It then tries every possible
color. Colors not yet used are equivalent under permutation, so only one new
color needs to be tried. A branch is discarded as soon as an uncolored
vertex has neighbors in all four colors.

\begin{algorithm}[H]
\caption{Exact four color search}
\label{alg:exact}
\begin{algorithmic}[1]
\Require Finite vertex set $V$ and forbidden pair set $E$
\Function{Search}{$c$}
    \If{every vertex is colored}
        \State \Return \textsc{true}
    \EndIf
    \State Choose an uncolored vertex $v$ of maximum saturation degree
    \State Break ties by the number of uncolored neighbors
    \State Let $T$ be the colors not used on colored neighbors of $v$
    \State Remove all but one color from $T$ that has not yet been used
    \For{each $q\in T$}
        \State Give $v$ color $q$
        \If{no uncolored vertex sees all four colors}
            \If{\Call{Search}{$c$}}
                \State \Return \textsc{true}
            \EndIf
        \EndIf
        \State Uncolor $v$
    \EndFor
    \State \Return \textsc{false}
\EndFunction
\State \Return \Call{Search}{all vertices uncolored}
\end{algorithmic}
\end{algorithm}

\FloatBarrier
\section{Conclusion}

We determined nine exact four color polynomial van der Waerden numbers and established
$\PVDW{ex^2+fx}{4}=\Theta(|e|+|f|)$
for every nonzero polynomial $ex^2+fx$. Our exact values are substantially smaller than the estimates given by general upper bounds, motivating the search for sharper bounds.

One immediate open problem is to determine the optimal constant $C_4$ such that
\[
\PVDW{ex^2+fx}{4}\le C_4(|e|+|f|)+1.
\]
We established that $57\le C_4\le64$.

More broadly, our work contributes to a research program of computing more exact polynomial van der Waerden numbers to better understand their growth. We aim to extend these methods to more colors and higher-degree polynomials, establish asymptotically tight bounds, and narrow the enormous gap between known lower and upper bounds. In particular, further exact values may help determine whether the growth of these numbers is closer to the exponential lower bounds or the much larger general upper bounds.

\end{document}